\documentclass[12pt]{article}
\usepackage{amsgen,amsmath,amstext,amsbsy,amsopn,amsfonts,amssymb}
\newtheorem{theorem}{Theorem}
\newtheorem{lemma}[theorem]{Lemma}
\newtheorem{corollary}[theorem]{Corollary}
\newtheorem{proposition}[theorem]{Proposition}
 \newtheorem{defi}[theorem]{Definition}

\newtheorem{exa}[theorem]{Example}

\newtheorem{rem}[theorem]{Remark}
\newenvironment{remark}{\begin{rem}\rm}{\end{rem}}
\newtheorem{rems}[theorem]{Remarks}

\def\bsq{\blacksquare\medskip}

\def\H{\mathcal H}
\def\L{\mathcal L}

\def\A{\mathcal A}

\def\X{\mathcal X}

\def\Hpi{(\pi, {\mathcal H})}
\def\unirep{unitary representation}
\def\supp{{\rm supp}}
\def\NN{{\mathbb N}}
\def\ZZ{{\mathbb Z}}
\def\CCC{{\mathbb C}}
\def\RRR{{\mathbb R}}

\def\RR+{{\mathbb R}^*}

\def\KK{{\mathbb K}}

\def\GG{{\mathbb G}}

\def\Q_p{{\mathbb Q}_p}
\def\S1{{\mathbb S}^1}

\def\Tr{{\rm Trace}}
\def\PT{Property (T)}
\def\eps{\varepsilon}
\def\Ga{\Gamma}

\def\vfi{\varphi}

\def\HS{{L^2}}
\def\spec{{\rm spec}}

\begin{document}

\title{A spectral gap property for random walks under unitary representations}
\author{ Bachir Bekka and Yves Guivarc'h}
\date{\today }

\maketitle

\begin{abstract}
Let $G$ be a locally compact group
and $\mu$ a probability measure on $G,$
which is not assumed to be  absolutely continuous
with respect to Haar measure.
Given a unitary representation
$\Hpi$ of $G,$ we study spectral properties
of the operator $\pi(\mu)$ acting
on $\H.$ Assume that $\mu$
is adapted and that the trivial representation
$1_G$ is not weakly contained in the tensor product
$\pi\otimes \overline\pi.$ We show that 
$\pi(\mu)$ has a spectral gap,
that  is, for the spectral radius $r_\spec(\pi(\mu))$ of 
$\pi(\mu),$ we have $r_\spec(\pi(\mu))<1.$
This provides a common generalization of several previously known
results.
Another consequence  is that, if $G$ has Kazhdan's Property (T), then
$r_\spec(\pi(\mu))<1$ for every unitary representation
$\pi$ of  $G$ without finite dimensional subrepresentations.
Moreover, we  give 
new examples of so-called identity excluding groups.
\end{abstract} 

\section{Introduction}
Let $G$ be a locally compact group
and $\mu$ a probability measure
defined on the Borel subsets of $G.$
We will always assume that $\mu$ is \emph{adapted};
 by this, we mean that the subgroup generated
by the support $\supp (\mu)$ of $\mu$ is dense
in $G.$ 
We will also consider
the stronger condition
that $\supp (\mu)$ is not contained
in the coset of a proper closed subgroup
of $G.$ In this case, we say that $\mu$ is \emph{strongly adapted}.
%(some authors say that $\mu$ is strictly aperiodic).

Let $\Hpi$ be a unitary representation
of $G,$ that is, a strongly continuous
homomorphism $\pi$ from $G$ to the unitary
group of a Hilbert space $\H.$
A bounded operator $\pi(\mu)$ on
$\H$ is defined by
$$
\langle\pi(\mu)\xi,\eta\rangle= \int_G \langle\pi(x)\xi,\eta\rangle d\mu(x), \qquad \forall \xi,\eta\in \H.
$$
Let $r_{\spec}(\pi(\mu))$ be the spectral radius of 
$\pi(\mu).$ We have $r_\spec(\pi(\mu))\leq 1,$
since  $\pi(\mu)$ is a contraction.

The main result of this paper is the
following theorem.
Recall that the (inner) tensor product
 $\pi\otimes \overline\pi$ is the representation acting on 
the Hilbert space $\H\otimes \overline\H$ by 
$\pi\otimes \overline\pi (x)=\pi(x)\otimes \pi (x)$ 
for all $x\in G.$
It can be realized
on the space $L^2(\H)$ of all Hilbert-Schmidt
operators on $\H$ by 
$$
\pi\otimes \overline\pi (x) T=\pi(x)T\pi(x^{-1}),\qquad x\in G,\, T\in \HS(\H).
$$
As is well-known and easy to see, 
$\pi\otimes \overline\pi$
has a non-zero invariant vector if
and only if $\pi$ contains a non-zero
finite dimensional subreprentation.

Recall also that a unitary representation
$\Hpi$ of $G$ almost has invariant vectors
if, for every compact subset $K$ of $G$
and every $\eps>0,$ there exists
a unit vector $\xi\in \H$ such that
$$
\max_{x\in K}\Vert \pi(x)\xi-\xi\Vert <\eps.
$$
Observe that  $\pi$ almost has invariant vectors
if and only if   the trivial
one-dimensional representation $1_G$
is weakly contained in $\pi,$
in the sense of \cite{Fell}.
 
\begin{theorem}
\label{Theo1} 
Let $G$ be a locally compact group
and $\mu$ an adapted probability measure
on $G.$ Let $\Hpi$ be a unitary representation
of $G.$ 
Assume that the following condition $(*)$
is satified:

$(*)$\qquad $\pi\otimes \overline\pi$ does not almost have 
invariant vectors.

\noindent
Then  $r_\spec(\pi(\mu))<1.$
If $\mu$ is strongly adapted, 
then  $\Vert\pi(\mu)\Vert <1.$ 
\end{theorem}

\begin{remark}
\label{Rem1}

(i) 
Observe that we 
do \emph{not} assume that $\mu$ is absolutely continuous with respect to 
a Haar measure on $G.$
Indeed, in this case, 
the result is well-known and easy to prove.
In fact, a stronger result is true:
if $\mu$ is absolutely continuous,
and strongly adapted
and if, instead of Condition $(*),$
we assume that
$\pi$ does not almost have 
invariant vectors, then 
$\Vert\pi(\mu)\Vert<1$
(see \cite[Proposition 4.1]{Guiv} or \cite[Appendix G.4]{BHV}).

(ii) In general,  we cannot replace Condition $(*)$
by the weaker condition that $\pi$ does not 
almost have invariant vectors.
For a simple counterexample, take
 $G=\S1$ the circle group, $\pi=\chi$ a non-trivial
character of $G,$ and  $\mu=\delta_z$ the Dirac measure
at a point $z$ generating a dense
subgroup of $\S1.$ In this case,
$r_\spec(\pi(\mu))=1.$

(iii)  In view of (ii),
one can ask, whether 
we have $r_\spec(\pi(\mu))<1,$
if $\pi$ does not 
almost have invariant vectors
and if we assume moreover
that $\mu$ is strongly adapted.
The answer is still negative as the
following example shows.
 Let $G=\S1$
and let $\pi$ be the restriction of
the regular representation of $\S1$ to
the subspace $L^2_0(\S1)$ 
of all functions in $L^2(\S1)$
which are orthogonal to the constants.
Let $\mu$ be any probability measure on $\S1$
with finite support. Then there exists
a sequence $(n_k)_k$ in $\ZZ\setminus\{0\}$
such that  $\lim_k z^{n_k}=1$ for all
$z\in \supp (\mu).$ This implies that
$1$ is a spectral value of  $\pi(\mu).$

(iv) Condition $(*)$ in the theorem above was studied 
by the first author in \cite{RepMoy}.
More specifically, a \unirep\ $\pi$
is  amenable in the sense of this paper 
if $\pi\otimes \overline\pi$ almost has invariant vectors.
In other words, the theorem can be stated as follows:
given an adapted probability measure $\mu$ on $G,$
a \unirep\ $\pi$ of  $G$  is amenable  if 
$r_\spec(\pi(\mu))=1.$

(v) Condition $(*)$  \emph{never} holds if $G$ is amenable.  
More precisely, it was shown
in \cite[Theorem 1]{BekkaTens} that, if $\pi$ is any unitary representation
of an amenable locally compact group, then $\pi\otimes \overline\pi$ almost has invariant vectors.

(vi) Actually, as the proof below shows,
we also obtain (under the assumptions of the theorem) the 
result that $r_\spec(\pi\otimes \overline\pi(\mu))<1$
(and  $\Vert\pi\otimes \overline\pi(\mu)\Vert<1$ in case $\mu$
is strongly adapted).

\end{remark}
\label{Rem5}
The theorem above leads to spectral gap properties
(see Corollaries~\ref{Cor2} and \ref{Cor3}) 
which are somewhat unexpected for general
probability measures.
This may be seen as an indication of the prevalence of exponential
mixing for actions of non amenable groups
(compare  \cite{Dolgopyat}).

We will derive several corollaries of our result.
The first one deals wih Kazhdan groups.
Recall that a locally compact group 
$G$ is said to have Kazhdan's \PT,
if whenever a unitary representation
$\pi$ of $G$ almost has invariant vectors,
then $\pi$ has a non-zero invariant vector
(see \cite{Kazhdan}, \cite{HaVa}, or \cite{BHV}). 
Examples of locally compact groups with
\PT\ are all simple real Lie  of real rank greater
or equal to $2.$  This is for instance the case for
the groups $SL_n(\RRR)$ for
$n\geq 3$ and $Sp_{2n}(\RRR)$ for $n\geq 2.$
An important fact is that
\PT\ is inherited by lattices.
(Recall that a discrete group $\Ga$ is
a lattice in a locally compact group $G$
if the homogeneous space $G/\Ga$ has a
$G$-invariant probability measure.)
Thus, $SL_n(\ZZ)$ for
$n\geq 3$ and $Sp_{2n}(\ZZ)$ for $n\geq 2$
have \PT. Examples of Lie groups
which have \PT\ and are not semisimple
are the semi-direct products 
$SL_n(\RRR)\ltimes \RRR^n$ for
$n\geq 3.$

The following corollary is an immediate consequence of the theorem
above.  
\begin{corollary}
\label{Cor2}
Let $G$ be a locally compact group with
\PT, and let $\mu$ be an adapted
probability measure on $G.$
Let $\pi$ be unitary representation
of $G$ without finite dimensional subrepresentations.
Then $r_\spec(\pi(\mu))<1.$
If $\mu$ is strongly adapted,
 then 
  $\Vert\pi(\mu)\Vert <1$.
\end{corollary}

\begin{remark}
\label{Rem4}
(i) Unfortunately, the previous corollary does not 
apply to the  important case
of a compact group; indeed,   every unitary representation of 
such a group is a direct sum of finite dimensional ones. 
It should be mentioned that, for the orthogonal groups
$SO(n,\RRR),\ n\geq 3,$ the
problem of finding a probability measure $\mu$
with a spectral gap for the corresponding  operator
acting on the subspace
orthogonal to the constants in $L^2({\mathbb S}^{n-1})$ is
 related to the so-called
Banach-Ruziewicz problem (see \cite{Lub}, \cite{Sar}).

(ii) The previous corollary  may be reformulated 
as follows, at least
when $G$ is second countable. Let $D$ be a dense subgroup 
of the Kazhdan group $G.$ For
any unitary representation $\pi$ of $G$ without finite
dimensional subrepresentation, the restriction
of $\pi$ to $D$ does not weakly contain
the trivial representation $1_D,$ where
$D$ is equipped with the discrete topology
(for more details, see Section~\ref{Section3}).
\end{remark}

Corollary~\ref{Cor2} can be generalized
to the case where $G$ has a closed subgroup $H$ such that
the pair $(G,H)$ has \PT.
Recall that $(G,H)$ has \PT\ (or, as some authors
say, $G$ has the relative \PT\ with respect to $H$)
if, whenever a unitary representation
$\pi$ of $G$ almost has invariant vectors,
then $\pi$ has a non-zero $H$-invariant vector
(see \cite{HaVa} or \cite{BHV}).
Observe that $G$ has \PT\ if and only
if the pair $(G,G)$ has \PT.
Another example of a pair with
\PT\ is the pair $(SL_2(\RRR)\ltimes \RRR^2, \RRR^2),$
for the semi-direct product $SL_2(\RRR)\ltimes \RRR^2$
with respect to the natural action of 
$SL_2(\RRR)$ on $\RRR^2.$

\begin{corollary}
\label{Cor5}
Let $G$ be a locally compact group
and $H$ a closed subgroup such that
the pair $(G,H)$ has
\PT, and let $\mu$ be an adapted
probability measure on $G.$
Let $\pi$ be unitary representation
of $G.$ Assume that the restriction of $\pi$ to $H$
has no finite dimensional subrepresentations.
Then $r_\spec(\pi(\mu))<1.$
If $\mu$ is strongly adapted,
 then 
  $\Vert\pi(\mu)\Vert <1$.
\end{corollary}

As we now see, Theorem~\ref{Theo1} above unifies and generalizes several
 previously known results. 
When $\pi$ is the left regular representation
$\lambda_G$ of $G$ acting on $L^2(G),$ 
we recover the result \cite[Th\'eor\`eme]{DeGu} 
by Derrienic and the second author (for a similar result,
see also \cite{BeCh}).
The case where $G$ is discrete and $\mu$
is symmetric is due to Kesten \cite{Kesten}.
Observe that $\lambda_G(\mu)$ is the  operator
$\xi\mapsto \mu\ast \xi$ on $L^2(G)$ given by left convolution with
$\mu.$
%$$\mu\ast \xi (x) =\int_G \xi(y^{-1}x) d\mu(y),\qquad \xi\in L^2(G),\ x\in G.$$

\begin{corollary}\textbf{(\cite{DeGu})} 
\label{Cor1}
Let $G$ be a non-amenable locally compact group,
and let $\mu$ be an adapted probability measure
on $G.$ Then $r_\spec(\lambda_G(\mu))<1.$
If $\mu$ is strongly adapted,
 then 
  $\Vert\lambda_G(\mu)\Vert <1$.
\end{corollary}

Our  proof of Theorem~\ref{Theo1} 
may be seen as a non-commutative
version of the  proof of \cite[Th\'eor\`eme]{DeGu}.
It relies on the characterization of amenable
representations given in \cite{RepMoy}
in terms of existence of invariant means on 
appropriate operator algebras
(see Theorem~\ref{Theo2} below).

The previous corollary
 generalizes to  homogeneous spaces
which are not amenable in the
sense of Eymard \cite{Eymard}. 
Even more generally, let $(X,\nu)$ be a 
measure space with a $\sigma$-finite measure
$\nu.$ Assume that 
the separable locally compact group
$G$ acts measurably on $X$ and 
that $\nu$ is quasi-invariant under
the action of $G.$
We say that the action of $G$ on $X$ 
is \emph{co-amenable} if 
there exists a $G$-invariant mean 
on $L^\infty(X,\nu),$ that is,
a positive linear functional $M$
on  $L^\infty(X,\nu)$ such that
that $M(1_X)=1$ and $M({}_g\vfi) =M(\vfi)$
for all $g$ in $G$ and $\vfi$ in $L^\infty(X,\nu),$
where ${}_g\vfi(x)=\vfi (g^{-1}x).$
(This notion should not be
confused with the well-established notion of 
an amenable action of $G$ on $X$
due to Zimmer;  see  \cite{Zimmer}.)
A unitary representation $\lambda_X$ of $G$
is defined on  $L^2(X,\nu)$ by 
$$\lambda_{X}(g) \xi(x)=\sqrt{\dfrac{d\nu(g^{-1}x)}{d\nu(x)}}\xi(g^{-1}x),
\qquad g\in G,\ x\in X,
\ \xi\in L^2(X).$$
The following result was obtained in
\cite[Proposition~1]{GuivAsym}
(see also \cite[Theorem 4.15]{Guiv}).

\begin{corollary}\textbf{(\cite{GuivAsym})} 
\label{Cor6}
Assume that the action of $G$
on $X$ is not co-amenable
and let $\mu$ be an adapted probability measure
on $G.$ Then $r_\spec(\lambda_X(\mu))<1.$
If $\mu$ is strongly adapted,
 then 
  $\Vert\lambda_X(\mu)\Vert <1$.
\end{corollary}

The previous corollary 
is a direct consequence of Theorem~\ref{Theo1}
and of the following proposition
which shows  that Condition $(*)$ 
is equivalent to the
non co-amenability of the action.

\begin{proposition} 
\label{Prop}
Let $G$ be   separable locally compact group
acting measurably on a measure space $(X,\nu),$ 
where $\nu$ is a $\sigma$-finite 
measure on $X$ which is quasi-invariant under
the action of $G.$
The following properties are equivalent:
\begin{itemize}
\item [(i)] the action of $G$ on $X$ is co-amenable;
\item [(ii)] the representation $\lambda_X $ almost has
invariant vectors;
\item [(iii)] the representation $\lambda_X \otimes \overline\lambda_X$ almost has
invariant vectors.
\end{itemize}
\end{proposition} 

The equivalence of (i) and (ii) 
in the previous proposition was also 
shown in \cite[Theorem 4.15]{Guiv}.

 \begin{remark}
\label{Rem2}
Assume that the action of  $G$ on the measure space $(X,\nu)$ 
is amenable in Zimmer's sense. The following remarkable
result was shown in \cite[Corollary 3.2.2]{Claire}:
$r_\spec(\lambda_X(\mu))=r_\spec(\lambda_G(\mu))$ for
any adapted probability measure $\mu$ on $G$
(the case of a discrete group was previously treated
in \cite{Kuhn}). In particular,
if $G$ is non amenable, then $r_\spec(\lambda_X(\mu))<1,$
by Corollary~\ref{Cor1}.
\end{remark}

We now turn to semisimple Lie groups.
In this case, using \cite[Lemma 4]{BekkaInvMean},
we immediately obtain the following strengthening of
Theorem~\ref{Theo1}. This result was shown
in \cite[Theorem~C]{ShalomFourier}.
The case of a semisimple Lie group  
with \PT\ was  treated in \cite[Theorem~1]{Nevo}.

\begin{corollary}\textbf{(\cite{ShalomFourier})} 
\label{Cor3}
Let $G$ be a semisimple real Lie group with
finite centre and without compact factors.
Let $\mu$ be an adapted
probability measure on $G.$
Let $\pi$ be unitary representation
of $G$  which does not almost have 
invariant vectors.
Then $r_\spec(\pi(\mu))<1.$
If $\mu$ is strongly adapted,
 then 
  $\Vert\pi(\mu)\Vert <1$.
\end{corollary}

Actually, the previous result
 is proved in \cite{ShalomFourier}
under the  weaker assumption that 
the projection of $\mu$ to every simple
factor of $G$  is not supported on a closed
amenable subgroup. It is worth mentioning
that the articles \cite{FurSha} and  \cite{Shalom}
contain interesting
 applications of the condition
$\Vert \pi(\mu)\Vert<1$ 
(or $\Vert \pi\otimes\overline\pi(\mu)\Vert<1$)
 to random ergodic theorems.

 Let $G$ be a separable locally compact
group acting measurably on a measure space
$(X,\nu),$ where $\nu$ is now
an  invariant probability measure.
As above, let $\lambda_X$ be the natural representation
of $G$ on $L^2(X).$ Observe that
$L^2(X)$ contains the constants functions
$\CCC 1_X.$
Let $\lambda_X^0$ denote the restriction
of $\lambda_X$  to the closed invariant subspace
$$
L^2_0(X)=\{\xi\in L^2(X)\ :\ \int_G f(x) d\nu(x) =0\}=(\CCC 1_X)^\perp.
$$
As is well-known, the action of $G$ is ergodic
if and only if $\lambda_X^0$ has no non-zero invariant vectors,
and
the action  is mixing if and only if $\lambda_X^0 \otimes \overline\lambda_X^0$ has no non-zero invariant vectors.
%Given an adapted measure $\mu$ on $G,$
%an important question is when does $r_\spec(\lambda_X^0(\mu)<1$
%(or $\Vert \lambda_X^0(\mu)\Vert <1$) hold?
Important examples of probability  spaces
$X$ as above are the spaces $X=G/\Gamma,$ 
where $\Gamma$ is a lattice in a locally compact group $G$
with action of $G$ by left translations.
When $G$
is a semisimple real Lie group, it was shown in
 Lemma~3 of \cite{BekkaInvMean} 
that, for any lattice $\Ga,$ the
representation $\lambda_X^0$
does not almost have invariant vectors.
As a consequence of this fact,  we obtain 
from Corollary \ref{Cor3} the following result,
which was also observed in \cite{ShalomFourier} 
and \cite{FurSha}.
 
\begin{corollary}
 \label{Cor4} Let $G$
be a semisimple real Lie group with
finite centre and without compact factors.
Let $\mu$ be  an adapted
probability measure on $G.$
Let $\Gamma$ be a lattice in $G$ 
and denote by $\lambda_{G/\Ga}^0$ the representation of $G$
on $L^2_0(G/\Gamma)$
Then $r_\spec(\lambda_{G/\Ga}^0(\mu))<1.$
\end{corollary}

As  a further appplication
of our results, we now show 
that some classes of groups have the so-called
 identity excluding property.
This property  
appeared in several articles
(see, for instance, \cite{JRT}, \cite{LinWitt},
\cite{JRW})
in connection with the study of the
behaviour of convolution powers of a
probability measure. 
It also plays an  important role in the study
of equidistribution properties of 
random walks (see \cite{GuivEqui}).
 It is defined as
follows. A locally compact group $G$
is \emph{identity excluding} 
if, for every irreducible
unitary representation $\pi$ of $G$ with
$\pi\neq 1_G$ and every  dense
subgroup $D$ of $G,$ the restriction
of $\pi$ to $D$ does not almost have
invariant vectors ($D$ being equipped with the
discrete topology).

Part (i) and Part (ii) of the following result are  direct consequences
of Corollary~\ref{Cor2} and  Corollary~\ref{Cor3}, and
Part (iii) requires a further argument (see Section~\ref{Section3}).

\begin{corollary}
 \label{Cor7} 
A secound countable locally compact group $G$
is identity excluding in the following cases:
\begin{itemize}
\item [(i)] $G$ has \PT;
\item [(ii)] $G$ is a semisimple real Lie group with
finite centre and without compact factors;
\item [(iii)] $G=\GG ({{\KK}})$ is the group of ${{\KK}}$-rational 
points of a ${{\KK}}$-isotropic simple algebraic group $\GG$ over 
a local field ${{\KK}}.$
\end{itemize}
\end{corollary}

This paper is organized as follows
In Section~\ref{Section1}, we give the
proof of Theorem~\ref{Theo1},
 in Section~\ref{Section2} the proof
of Proposition~\ref{Prop}, and in
Section~\ref{Section3} the proof
of Corollary~\ref{Cor7}.
\section{Proof of Theorem~\ref{Theo1}}
\label{Section1}
Our proof of Theorem~\ref{Theo1} is modelled
after the proof of \cite[Th\'eor\`eme]{DeGu}.
Let $\H$ be a Hilbert space and $\L(\H)$  the algebra
 of all bounded operators
on $\H.$ Let $\A$ be a $C^*$-subalgebra
of $\L(\H),$ that is, a norm closed self-adjoint
subalgebra of $\L(\H).$
Assume that $\A$ is unital, that is, it contains 
the identity operator $I.$
A state   $M$ on $\A$ is
a  linear functional on $\A$ which is positive
(that is, $M(T^*T)\geq 0$ for all $T\in \A$)
and normalized (that is, $M(I)=1$).
In view of the analogy with
the classical case, one should
think of $M$ as a mean on $\A.$ 

Let now  $\pi$ be a unitary representation
of the locally compact group $G$ on $\H.$
We denote by ${\rm Ad}$ the representation
of $G$ on the vector space $\L(\H)$ given by
$$
{\rm Ad}(x)(T)=\pi(x) T\pi(x^{-1})\qquad x\in G,\ T\in \L(\H).
$$
Observe that this representation is not necessarily strongly continuous,
when $\L(\H)$ is endowed with the norm topology; this means
that, for fixed
$T\in\L(\H),$
the mapping 
$$
G\to \L(\H),  \quad x\mapsto \pi(x) T\pi(x^{-1})
$$
is in general not continuous (as an easy example, take $G=\RRR,$  $\pi=\lambda_{\RRR}$ the regular representation
of $L^2(\RRR),$ and $T$ the mutiplication operator
 defined by a bounded
function  which is not uniformly continuous). 

Let $\X(\H)$ be the subspace of $\L(\H)$ consisting
of all $G$-continuous operators, that is, all
$T\in \L(\H)$ such that the mapping 
$x\mapsto \pi(x) T\pi(x^{-1})$ is continuous.
It is easy to show that $\X(\H)$ is an 
${\rm Ad}(G)$-invariant unital $C^*$-subalgebra
of $\L(\H)$. Using Cohen's factorization theorem
(see, e.g., \cite[Theorem~(32.22)]{HewRos}
we have the following description of 
$\X(\H):$
$$
\X(\H)=\{{\rm Ad}(f)T\ :\ f\in L^1(G),\ T\in \L(\H)\},
$$
where ${\rm Ad}(f)T$ is defined by 
means of the weakly convergent integral
$$
{\rm Ad}(f)T= \int_G f(x) \pi(x)T\pi(x^{-1}) dx.
$$
Let $\A$ be an  ${\rm Ad}(G)$-invariant unital $C^*$-subalgebra
of $\L(\H).$ A state $M$ on $\A$ is ${\rm Ad}(G)$-invariant 
if
$$
M(\pi(x)  T\pi(x^{-1})=M(T)\qquad \forall\  T\in \A,\ x\in G.
$$

The essential tool for our proof of Theorem~\ref{Theo1}
is the following result from \cite{RepMoy}
which is a generalization
of the Hulanicki-Reiter characterization
of amenable groups (see \cite[3.2.5]{Green}). 
\begin{theorem}(\textbf{\cite{RepMoy}})
\label{Theo2}
 Let   $\pi$ be a unitary representation
of the locally compact group $G$ on $\H.$
The following properties are equivalent:
\begin{itemize}
\item[(i)]  $\pi\otimes\overline\pi$ almost has invariant vectors;
\item[(ii)] there exists an ${\rm Ad}(G)$-invariant 
state on $\L(\H);$
\item[(iii)] there exists an ${\rm Ad}(G)$-invariant 
state on the algebra $\X(\H)$ of all $G$-continuous operators.
\end{itemize}
\end{theorem}

For the proof of Theorem~\ref{Theo1},
we need two elementary lemmas. The first one is well-known 
and easy to prove (see the beginning of the proof
of Th\'eor\`eme~1 in \cite{DeGu}). 
\begin{lemma}
\label{Lem1}
Let $T$ be a linear contraction on a Hilbert space
$\H.$ Let $c\in \CCC$ with $\vert c\vert =1$
be a spectral value of $T.$ Then $c$ is an approximate
eigenvalue of $T,$
that is, there exists a sequence $(\xi_n)_n$ of unit vectors
 in $\H$ such that $\lim_n\Vert T\xi_n-c\xi_n\Vert=0.$
\end{lemma}

For a Hilbert space
$\H,$ recall that the inner product
on the space  $\HS(\H)$ of all Hilbert-Schmidt operators
on $\H$ is defined by 
$$
\langle S,T\rangle =\Tr (T^*S), \qquad S, T\in \HS(\H),
$$
where $\Tr$ is the usual trace of operators.
We denote by $T\mapsto \Vert T\Vert_2$ the corresponding norm.
Recall that the absolute value of an operator 
$T\in \L(\H)$ is the positive
operator $\vert T\vert =(T^*T)^{1/2}.$
\begin{lemma}
\label{Lem2} Let $S, T\in \HS(\H).$ Then
$$
\vert\langle S,T\rangle\vert^2 \leq \langle \vert S\vert,\vert T\vert \rangle
\langle \vert S^*\vert,\vert T^*\vert \rangle.
$$
\end{lemma}
\begin{proof}
Let $U$ and $V$ be the partial isometries
on $\H$ arising from the
polar decompositions 
$$S=U\vert S\vert, \qquad T=V\vert T\vert.$$
%Observe that $\vert\Tr (AB)\vert \leq \Vert A\Vert  \Tr (\vert B\vert),$for any bounded %operator $A$ on $\H$ and any  trace-class operator $B,$ where $\Vert A\Vert$ denotes
%the operator norm of $A.$
Observe that the polar decompositions of $S^*$ and $T^*$ are
$$
S^*=U^*( U\vert S\vert U^*), \qquad T=V^*( V\vert T\vert V^*),
$$
so that $\vert S^*\vert =U\vert S\vert U^*$ and 
$\vert T^*\vert =V\vert T\vert V^*.$
Using the Cauchy-Schwarz inequality
in $\HS(\H),$ we have
\begin{eqnarray*}
\vert\langle S,T\rangle\vert &=&
\vert\Tr (V^*\vert T^*\vert U \vert S\vert)\vert\\
&=&
\vert\Tr (V^*\vert T^*\vert^{1/2}\vert T^*\vert^{1/2} U 
\vert S\vert^{1/2}\vert S\vert^{1/2})\vert\\
&=&
\vert\Tr \left((\vert S\vert^{1/2}V^*\vert T^*\vert^{1/2})
(\vert T^*\vert^{1/2} U \vert S\vert^{1/2})\right)\vert\\
&\leq&
\Vert \vert S\vert^{1/2}V^*\vert T^*\vert^{1/2}\Vert_2 \ 
\Vert \vert T^*\vert^{1/2} U \vert S\vert^{1/2}\Vert_2\\
&=&
\left(\Tr (\vert T^*\vert^{1/2} V \vert S\vert V^*\vert T^*\vert^{1/2})\right)^{1/2}
\left(\Tr (\vert S\vert^{1/2} U^* \vert T^*\vert U\vert S\vert^{1/2})\right)^{1/2}\\
&=&
\left(\Tr (\vert V^*\vert T^*\vert V \vert S)\right)^{1/2}
\left(\Tr (\vert U\vert S\vert U^* \vert T^*)\right)^{1/2}\\
&=&
\left(\Tr (\vert T\vert \vert S\vert)\right)^{1/2}
\left(\Tr (\vert S^*\vert  \vert T^*\vert )\right)^{1/2}\\
&=&
\left(\langle \vert S\vert,\vert T\vert \rangle\right)^{1/2}
\left(\langle \vert S^*\vert,\vert T^*\vert \rangle\right)^{1/2}.
\end{eqnarray*}
$\bsq$
\end{proof}

We now proceed with the proof of Theorem~\ref{Theo1}.
This will be done in several steps.

\noindent
$\bullet$ {\it First step:}
We claim that it suffices to show  that $r_\spec((\pi\otimes \overline\pi)(\mu))<1.$
To show this, we use the following  argument which appears in
the proof of Theorem~1 in \cite{Nevo}. For every   $n\in \NN$ and  every vector $\xi\in\H,$ 
denoting by $\mu^{*n}$ the $n$-fold convolution product
of $\mu$ by itself, we have
\begin{eqnarray*}
\vert\langle(\pi(\mu))^{n}\xi,\xi \rangle\vert^2 
&=&\vert\langle\pi(\mu^{*n})\xi,\xi \rangle\vert^2\\
&=&\left\vert\int_G \langle\pi(x)\xi,\xi\rangle
d\mu^{*n}(x)\right\vert^2\\
&\leq&\int_G \vert\langle\pi(x)\xi,\xi \rangle\vert^2
d\mu^{*n}(x)\\
&=&\int_G \langle\pi\otimes \overline\pi(x)\xi\otimes \overline\xi,\xi\otimes \overline\xi \rangle
d\mu^{*n}(x)\\
&=&\langle(\pi\otimes \overline\pi(\mu))^{n}\xi\otimes \overline\xi,\xi\otimes \overline\xi \rangle, 
\end{eqnarray*}
where we used Jensen's inequality.
It follows that
$$
\Vert\pi(\mu)^{n}\Vert \leq \Vert(\pi\otimes \overline\pi(\mu))^{n}\Vert^{1/2}
$$
for all $n\in\NN$ and hence
$$
r_\spec((\pi)(\mu))=\lim_n \Vert\pi(\mu)^{n}\Vert^{1/n}\leq
\left(r_\spec((\pi\otimes \overline\pi)(\mu))\right)^{1/2}.
$$
This proves the claim.

From now on, we assume that $r_\spec((\pi\otimes \overline\pi)(\mu))=1.$
We will show that $\pi\otimes \overline\pi$  almost has
invariant vectors.

$\bullet$ {\it Second step:}
We claim that $1$ is an approximate
eigenvalue of $(\pi\otimes \overline\pi)(\mu).$ 
Indeed, by assumption, the operator
 $(\pi\otimes \overline\pi)(\mu)$ on $\HS(\H)$ has a
spectral value $c$ with $|c|=1.$ 
 Since $(\pi\otimes \overline\pi)(\mu)$ is a contraction, 
there exists, by Lemma~\ref{Lem1},  a sequence of operators $T_n$ in $\HS(\H)$ 
with $\Vert T_n \Vert_2=1$ such that
$$
\lim_n \Vert  (\pi\otimes \overline\pi)(\mu)(T_n)-cT_n \Vert_2 =0, 
$$ 
that is,
$$
\lim_n \left\Vert \int_G\pi(x)T_n\pi(x^{-1})d\mu(x)-c T_n\right\Vert_2=0,
$$
or, equivalently,
$$
\lim_n\int_G \langle \pi(x)T_n\pi(x^{-1}),T_n\rangle d\mu(x) =
\lim_n \langle(\pi\otimes \overline\pi)(\mu)(T_n), T_n \rangle =c.
$$
Using Lemma~\ref{Lem2} and Cauchy-Schwarz inequality, we have
on the one hand
\begin{eqnarray*}
&&\left|\int_G \langle \pi(x) T_n\pi(x^{-1}),T_n \rangle d\mu(x)\right|
\leq\int_G  |\langle \pi(x) T_n\pi(x^{-1}),T_n \rangle| d\mu(x)\\
&\leq&\int_G  
\langle \pi(x) |T_n|\pi(x^{-1}),|T_n|\rangle^{1/2} \langle \pi(x) |T_n^*|\pi(x^{-1}),|T_n^*|\rangle^{1/2}d\mu(x)\\
&\leq&\left(\int_G  
\langle \pi(x) |T_n|\pi(x^{-1}),|T_n|\rangle d\mu(x)\right)^{1/2}
\left(\int_G  
\langle \pi(x) |T_n^*|\pi(x^{-1}),|T_n^*|\rangle d\mu(x)\right)^{1/2}.\\
\end{eqnarray*}
On the other hand, we have
$$
\int_G  \langle \pi(x) |T_n|\pi(x^{-1}),|T_n|\rangle d\mu(x) \leq 
\int_G  \Vert\pi(x)T_n\pi(x^{-1})\Vert_2 \Vert T_n\Vert_2 d\mu(x) =1
$$
as well as 
$$
\int_G  \langle \pi(x) |T_n^*|\pi(x^{-1}),|T_n^*|\rangle d\mu(x) \leq 
\int_G  \Vert\pi(x)T_n^*\pi(x^{-1})\Vert_2 \Vert T_n^*\Vert_2 d\mu(x) =1.
$$
Since
$$
\lim_n \left|\int_G \langle \pi(x) T_n\pi(x^{-1}),T_n \rangle d\mu(x)\right| =1,
$$
it follows that
$$
\lim_n \langle(\pi\otimes \overline\pi)(\mu)(\vert T_n\vert), \vert T_n\vert \rangle=
\lim_n\int_G  \langle \pi(x) \vert T_n\vert\pi(x^{-1}),\vert T_n\vert\rangle d\mu(x)= 1,
$$
that is, 
$$
\lim_n \Vert (\pi\otimes \overline\pi)(\mu)(\vert T_n\vert)-\vert T_n\vert\Vert_2=0.
$$
Since $\Vert \vert T_n\vert\Vert_2=1,$ this shows
that $(\pi\otimes \overline\pi)(\mu)$ 
has $1$ as approximate eigenvalue. 
This proves the second step.

$\bullet$ {\it Third step:}
There exists a state $M$ on 
$\L(\H)$ which is invariant under a
dense subgoup $D$ of $G.$
Indeed, by the second step, there exists
a sequence $(T_n)_n$ of  Hilbert-Schmidt operators
on $\H$ with $\Vert  T_n\Vert_2=1$ 
such that
$$
\lim_n\int_G  \langle \pi(x) T_n\pi(x^{-1}), T_n\rangle d\mu(x)= 1,
$$
It follows that there exists a subsequence, still denoted by 
$(T_n)_n,$ such that
$$
 \lim_n  \langle \pi(x) T_n\pi(x^{-1}), T_n\rangle=1
$$
and therefore 
$$
 \lim_n \Vert \pi(x) T_n\pi(x^{-1})-T_n\Vert_2= 0 \leqno{(*)}
$$
for $\mu$-almost every $x$ in $G.$ The set of all
$x$ for which $(*)$ holds is clearly a measurable subgroup $D$ of $G.$
Since $\mu(D)=1,$ the support of $\mu$ is contained in the
closure $\overline{D}$ of $D.$
By assumption,
the support of $\mu$ generates a dense
subgroup of $G.$ Hence, $D$ is dense in $G.$

Consider now the sequence of states $M_n$ on $\L(\H)$ defined by
$$
M_n(T)= \langle TT_n, T_n\rangle, \qquad T\in \L(\H).
$$
Let $M$ be a weak *-limit point of $(M_n)_n$ in the
dual space of $\L(\H).$
Using $(*)$ above, it is readily verified that $M$ is a state
on $\L(\H)$ which is ${\rm Ad}(D)$-invariant.

$\bullet$ {\it Fourth step:} the representation
$\pi\otimes\overline \pi$ almost has invariant vectors.
Indeed, consider the restriction $M_0$ of
$M$ to the subalgebra $\X(\H)$ of all $G$-continuous operators
on $\H.$ Then $M_0$ is an ${\rm Ad}(D)$-invariant
state on $\X(\H).$ Since $D$ is dense in $G,$
it follows from the continuity 
of the action of $G$ on $\X(\H)$
that $M_0$ is 
${\rm Ad}(G)$-invariant. Theorem~\ref{Theo2}
shows then that $\pi\otimes\overline \pi$ almost has invariant vectors.
This concludes the proof of 
the main part of Theorem~\ref{Theo1}.

To show the last part
of Theorem~\ref{Theo1},
observe first that 
 $\mu$ is strongly adapted
if and only if the probability measure
$\nu=\check{\mu}\ast\mu$ is adapted,
 where
the  measure $\check\mu$ is defined by 
$$\check \mu(E)={\mu(E^{-1})}$$
for every Borel subset $E$ of $G.$
Observe also  $\pi(\nu)=\pi(\mu)^*\pi(\mu)$ is a self-adjoint positive
operator with norm $1$ and that
$\Vert\pi(\nu)\Vert=\Vert\pi(\mu)\Vert^2.$

Assume now that $\Vert\pi(\mu)\Vert=1.$
Then $\Vert\pi(\nu)\Vert=1$ and hence
$1$ is a spectral value of $\pi(\nu).$
The proof above applied to 
 $\nu$ in place of $\mu$ shows that 
$\pi\otimes\overline \pi$ almost has invariant vectors.
$\bsq$

\section{Proof of Proposition~\ref{Prop}}
\label{Section2}
Let $G$ be a separable locally compact group
acting  on the measure space $X$ equipped with a quasi-invariant 
$\sigma$-finite measure
$\nu.$ Let $\lambda_X$  be the 
 corresponding unitary representation on
  $L^2(X,\nu).$
%Assume that $r_{\spec}(\lambda_X(\mu))=1$
%for an adapted probability measure $\mu$ on $G.$ 
%By Theorem~\ref{Theo1}, the representation $\lambda_X\otimes\overline\lambda_X$
%almost has invariant vectors.
%Assume that $r_{\spec}(\lambda_X(\mu))=1$
%for an adapted probability measure $\mu$ on $G.$ 
To show  that (i) implies (ii), assume that the action of
$G$ on $(X,\nu)$ is co-amenable.
Thus, there exists a $G$-invariant mean
on $L^\infty(X,\nu).$
Then, by standard arguments, there exists a sequence
$(f_n)_n$ of positive measurable functions
in $L^1(X)$ with $\Vert f_n\Vert_1=1$ such that
$\lim_n \Vert\lambda_X(g) f_n-f_n\Vert_1 =0$ uniformly
on compact subsets
of $G.$ Set $\xi_n=\sqrt{f_n}.$
Then $\Vert \xi_n\Vert_2=1$ and
$\lim_n \Vert\lambda_X(g) \xi_n-\xi_n\Vert_2 =0$ uniformly
on compact subsets
of $G,$ so that $\lambda_X$
almost has invariant vectors.

The fact that (ii) applies (iii) is obvious.
To show that (iii) implies (i),
assume that 
the representation $\lambda_X\otimes\overline\lambda_X$
almost has invariant vectors.
It follows from Theorem~\ref{Theo2}
that there exists an ${\rm Ad}(G)$-invariant state $M$
on the algebra $\L(L^2(X,\nu)).$
Observe that the algebra $L^\infty(X,\nu)$ can be embedded
as a subalgebra of $\L(L^2(X,\nu)):$
to every $\vfi\in L^\infty(X,\nu)$ corresponds
the operator $T_\vfi$ acting on $L^2(X,\nu)$ 
by pointwise mutiplication by $\vfi.$
For $g\in G$ and $\vfi\in L^\infty(X,\nu),$
we have the relation
$$
\lambda_X(g) T_\vfi\lambda_X(g^{-1})= T_{{}_g\vfi} \leqno{(**)}
$$
where ${}_g\vfi(x)=\vfi(g^{-1}x)$.

Define now a mean $M_0$ on $L^\infty(X,\nu)$
by
$$
M_0(\vfi)= M(T_\vfi),\qquad \forall \vfi\in L^\infty(X,\nu).
$$
Relation $(**)$ above shows that $M_0$
is $G$-invariant. Hence, the action of
$G$ on $(X,\nu)$ is co-amenable.

\section{Proof of Corollary~\ref{Cor7}}
\label{Section3}

Before we proceed with the proof of Corollary~\ref{Cor7},
several remarks are in order.
In the following lemma, we will use
and prove only the fact that (i) implies (ii).
The fact that (ii) implies (i),
which was  shown in \cite[Theorem~3.6]{JRT}
(see also \cite[Proposition~1]{GuivAsym}), is 
actually the main motivation for the
notion of an identity excluding group.
The proof of this implication 
uses the arguments which appeared
in the Third Step of the proof of Theorem~\ref{Theo1}.
\begin{lemma}
\label{Lem3}
Let $G$ be a second countable, locally compact group.
Let $\Hpi$ be a unitary representation of $G.$
The following properties are equivalent:
\begin{itemize}
\item[(i)] There exists a dense  subgroup $D$ of $G$
such that the restriction $\pi\vert_D$ of $\pi$ to $D$
almost has invariant vectors;
\item[(ii)] there exists an adapted 
 probability measure $\mu$ on the Borel subsets of  $G$
such that  $1$ is an approximate eigenvalue of $\pi(\mu).$
\end{itemize}
\end{lemma}

\begin{proof}
To show that (i) implies (ii), 
let $(U_n)_{n\in \NN}$ be  a basis
for the topology of $G.$ Since 
$D$ is dense, we can 
find, for every $n,$ and  element $g_n\in D\cap U_n.$
Clearly, the sequence $(g_n)_n$ is dense in $G.$
Set 
$$
\mu=\sum_{n\in \NN} 2^{-n} \delta_{g_n}.
$$
Then $\mu$ is  an adapted  probability measure on $G.$
Since $\pi\vert_D$ 
almost has invariant vectors, there exists
a sequence $(\xi)_i$ of unit vectors 
in $\H$ such that, for every $n\in \NN,$ 
$$
\lim_i \Vert \pi(g_n)\xi_i-\xi_i\Vert =0.
$$
It follows that $1$ is an approximate eigenvalue
of $\pi(\mu).$ Indeed, for a given $\eps>0,$
choose $N$ such that $\sum_{n>N} 2^{-n}\leq \eps.$
Let $i_0$ be such that, for all $n\leq N,$ 
$$
\Vert \pi(g_n)\xi_i-\xi_i\Vert \leq \eps/N \qquad \text{for all}\quad i\geq i_0.
$$
Then, for all $i\geq i_0,$ we have
\begin{eqnarray*}
 \Vert \pi(\mu)\xi_i -\xi_i\Vert
&\leq& \sum_{n\geq 1} \Vert \pi(g_n)\xi_i-\xi_i\Vert\\
&\leq& \sum_{n\leq N} \Vert \pi(g_n)\xi_i-\xi_i\Vert + 2\eps\\
&\leq& 3\eps.\quad\bsq
\end{eqnarray*}
\end{proof}

\begin{remark}
\label{Rem3}
 An irreducible unitary representation $\Hpi$ of  
a locally compact group $G$ is 
said to be CCR if, for every $f\in L^1(G),$
the operator $\pi(f)$ is a compact operator
on $\H.$
If this is the case,
 then $\pi$ is -up to unitary equivalence- the unique irreducible
unitary representation of $G$ which is weakly
contained in $\pi$ 
(see \cite[Corollaire 4.1.11]{Dixmier}). In particular, a CCR representation
$\pi\neq 1_G$ does not almost have invariant vectors.
This shows that Part (i) of Corollary~\ref{Cor7}
is a consequence of Corollary~\ref{Cor2}.

It is well-known that  all  irreducible unitary representations
of a semisimple real Lie group
with finite centre
are CCR (see \cite[Th\'eor\`eme 15.5.6]{Dixmier}).
\end{remark}

We now give the proof
of Corollary \ref{Cor7}. In view of the previous remark and of Lemma~\ref{Lem3},
it is clear that Part (i) and Part (ii) 
follow from Corollary \ref{Cor2} and Corollary \ref{Cor3},
respectively.
It remains to prove Part (iii).

So, let $G=\GG ({{\KK}})$ be the group of ${{\KK}}$-rational 
points of a ${{\KK}}$-isotropic simple algebraic group $\GG$ over 
a local field ${{\KK}}.$ Let $\pi$ be an irreducible unitary representation of
$G.$ We claim that  $r_\spec(\pi(\mu))<1$ for every adapted 
probability measure $\mu$ on $G.$

Indeed, it is known that 
there is a real number $p$ in $[2, \infty)$
such that all the matrix coefficients of $\pi$ lie in $L^p(G)$
(see \cite[Chap.~XI, 3.6 Proposition]{BW}).
Hence, for some integer $N$, the  tensor power
$\pi^{\otimes N}\otimes\overline {\pi}^{\otimes N}$ is contained in an infinite multiple
of the regular representation $\lambda_G$.
It follows that
$$
r_\spec((\pi^{\otimes N}\otimes\overline {\pi}^{\otimes N})(\mu))\leq r_\spec(\lambda_G(\mu)).
$$
On the other hand, the argument used in the First Step of the proof of 
Theorem~\ref{Theo1}
shows that
$$
r_\spec((\pi)(\mu))\leq 
\left(r_\spec((\pi^{\otimes N}\otimes \overline{\pi}^{\otimes N})(\mu))\right)^{1/2N}.
$$
Hence $ r_\spec((\pi)(\mu))\leq r_\spec(\lambda_G(\mu))^{1/2N}.$
As $G$ is not compact, $G$ is not amenable.
Therefore, $r_\spec(\lambda_G(\mu))<1,$ by Corollary~\ref{Cor1}.
This proves the claim.

\goodbreak
\noindent
{\bf Addresses}

\noindent
 Bachir Bekka, UFR Math\'ematique, Universit\'e de  Rennes 1, 
Campus Beaulieu, F-35042  Rennes Cedex, France

\noindent
E-mail : bachir.bekka@univ-rennes1.fr

\noindent
Yves Guivarc'h, UFR Math\'ematique, Universit\'e de  Rennes 1, 
Campus Beaulieu, F-35042  Rennes Cedex, France

\noindent
E-mail : yves.guivarch@univ-rennes1.fr

\end{document}